\documentclass[12pt]{article}
\usepackage[ansinew]{inputenc}
\usepackage[pdftex]{graphicx}
\usepackage{amssymb,amsmath,amscd}
\usepackage{dsfont}
\usepackage{verbatim}
\usepackage{amsmath} 
\usepackage{graphicx}
\usepackage{pstricks}
\usepackage{enumitem}
\usepackage{ntheorem}
\usepackage{amssymb}
\usepackage{float,caption}

\newtheorem{de}{Definition}[section] 
\newtheorem{nott}{Notation}
\newtheorem{ex}[de]{Example}
\newtheorem{theo}{Theorem}
\newtheorem{prop}[theo]{Proposition}  
\newtheorem{lem}{Lemma}[section]
\newtheorem{cor}[lem]{Corollary}

\setenumerate[1]{font=\bfseries,label=\Roman*.}
\setenumerate[2]{font=\itshape,label=\arabic*)}

\author{Yashar Memarian\\
        Laboratoire de Math\'ematiques d'Orsay,\\
        Univ Paris-Sud,\\
        Orsay, F-91405,\\
        CNRS,Orsay,France}

\title{A Lower Bound on the Waist of Unit Spheres of Uniformly Convex Normed Spaces\\
       }
\date{}

\begin{document}
\maketitle

\begin{abstract}
In this paper we give a lower bound on the waist of the unit sphere of a uniformly convex normed space by using the localization technique in 
codimension greater than one and a strong version of the Borsuk-Ulam theorem. The tools used in this paper follow ideas of M. Gromov in \cite{grwst}. Our isoperimetric type inequality generalizes the Gromov-Milman isoperimetric inequality in \cite{gromil}.
\end{abstract}

\section{introduction}

\footnotetext{E-mail address: yashar.memarian@math.u-psud.fr}

The classical isoperimetric inequality for a metric space relates the measure of compact sets to the measure of their boundaries. These inequalities are codimension $1$ isoperimeric inequalities (simply because the difference of the dimension of a compact set and the dimension of its boundary is equal to $1$).
\\ During his research on a Morse theory for the space of cycles of a manifold, F. Almgren gave a sharp lower bound for the volume of a minimal $k$-cycle in the sphere $\mathbb{S}^n$ for every $k$ (see \cite{pitz},\cite{grfil}). This is an instance of an higher codimensional isoperimetric type inequality.
\\ Another important example of higher codimensional isoperimetric inequality, which in fact is a generalisation of the Almgren isoperimetric inequality on the sphere, is the waist of the sphere theorem of Gromov presented in \cite{grwst}.
\\ In this paper we prove a higher codimensional isoperimetric inequality for the unit sphere of a uniformly convex normed space. The idea follows \cite{grwst}.
\\ In \cite{gromil}, M. Gromov and V. Milman give an isoperimetric inequality for the unit sphere of a uniformly convex normed space by using the localization technique (a nice exposition of this can be found in \cite{ale}). The main result of this paper, generalizes the isoperimetric inequality of Gromov-Milman.
\\ We begin by defining waist. For more details about this invariant see \cite{grwst},\cite{memwst}.
\begin{nott}[Tubular neighborhoods]
Let $X$ be a metric space, $Y$ a subset of $X$, $\varepsilon>0$. The $\varepsilon$-neighborhood of $Y$ is denoted by
$$Y+\varepsilon = \{x \in X \,\vert\, d(x,Y) \leq \varepsilon \}.$$
\end{nott}

\begin{de}[Waist of a metric-measure space, see \cite{grwst}]
Let $X=(X,d, \mu)$ be a mm-space. Let $Z$ be a topological space. Let $w(\varepsilon)$ be a positive function. We say the {\em waist of $X$ relative to $Z$} is larger than $w$ if for every continuous map $f : X \rightarrow Z$ there exists a $z \in Z$ such that for all $\varepsilon>0$,
$$\mu(f^{-1}(z)+ \varepsilon) \geq w(\varepsilon).$$
\end{de}

The purpose of this paper is to give a lower bound of the waist of the unit sphere of a uniformly convex normed space relative to $\mathbb{R}^k$. We are ready to state the main theorem of this paper.

\begin{theo} \label{usphere}
Let $X$ be a uniformly convex normed space of finite dimension $n+1$. Let $S(X)$ be the unit sphere of $X$, for which the distance is induced from the norm of $X$. The measure defined on $S(X)$ is the conical probability measure. Then a lower bound for the waist of $S(X)$ relative to $\mathbb{R}^{k}$ is given by
$$w(\varepsilon)= \frac{1}{1+(1-2\delta(\frac{\varepsilon}{2}))^{n-k}(k+1)^{k+1}\frac{F(k,\frac{\varepsilon}{2})}{G(k,\frac{\varepsilon}{2})}}$$
where $\delta(\varepsilon)$ is the modulus of convexity,
$$F(k,\varepsilon)= \int_{\psi_2(\varepsilon)}^{\frac{\pi}{2}}\sin(x)^{k-1}\,dx.$$ 
and 
$$G(k,\varepsilon)= \int_{0}^{\psi_1(\varepsilon)}\sin(x)^{k-1}\,dx.$$
And where 
$$\psi_1(\varepsilon)=2\arcsin(\frac{\varepsilon}{4\sqrt{k+1}})$$
and
$$\psi_2(\varepsilon)=2\arcsin(\frac{\varepsilon}{2\sqrt{k+1}})$$
\end{theo}

Section 2 will be concerned with preliminaries and tools which we need to prove this theorem. In the last section we will discuss the relation of our result with Gromov-Milman's isoperimetric inequality and some applications of our theorem.

\subsection{Acknowledgement}
We thank S. Alesker whose exposition \cite{ale} has helped us a lot.

\section{Preliminaries}

Let us consider a uniformly convex normed space of dimension $(n+1)$, $X=(\mathbb{R}^{n+1}, \Vert \quad \Vert)$ which we fix once for all.

\begin{de}[Modulus of convexity]
The space $X$ has modulus of convexity $\delta$ if for all $\varepsilon>0$, for all vectors $x,\,y \in\ X$ with $\Vert x \Vert=\Vert y \Vert=1$ and 
$\Vert x-y \Vert \geq \varepsilon$ we have 
$$\frac{\Vert x+y \Vert}{2} \leq 1-\delta(\varepsilon).$$
\end{de}

\begin{ex}
Let $E$ be a Euclidean space. In this case, the modulus of convexity is easily determined from the parallegram identity. And we have
$$\delta_{E}(\varepsilon)=1-\sqrt{1-\frac{\varepsilon^2}{4}}.$$
\end{ex}

\emph{Remark} $\delta$ is a monotone increasing function. We use this remark later on to prove the Lemma $5.2$. 

We denote by $B(X):=\{x \in\ X \vert \quad \Vert x \Vert \leq 1\}$ the unit ball of $X$ and $\partial{B(X)}=S(X):=\{x \in\ X \vert \quad \Vert x \Vert=1\}$
the unit sphere of $X$.
\\ We define a probability measure $\mu$ on $S(X)$ and we call it the conical measure,

\begin{de}[conical probability measure]
For any Borel set  $A \subset S(X)$ we define 
$$\mu(A):=\frac{m_{n+1}\{\bigcup tA \vert 0 \leq t \leq 1\}}{m_{n+1}(B(X))}$$
where $m_n$ is the $n$-dimensional Lebesgue measure on $X$.
\end{de}
We can check that the measure $\mu$ is a probability measure on $S(X)$, indeed
$$\mu(S(X))=\frac{m_{n+1}\{tS(X), 0 \leq t \leq 1\}}{m_{n+1}B(X)}=1.$$

\emph{Remark}: For the Euclidean norm on $\mathbb{R}^{n+1}$, where the distance between two points is the Euclidean 
distance and where the unit sphere is the canonical $n$-dimensional sphere $\mathbb{S}^n$, the conical measure is the 
canonical Riemannian probability measure on $\mathbb{S}^n$.

The mm-space on which we are going to work is 
$(S(X),\mu, d)$ with $\mu$ the conical probability measure and $d$ the distance induced on $S(X)$
from the norm defined on $X$ (i.e. for all $x,y \in S(X)$, $d(x,y)=\Vert x-y \Vert$).

\section{Scheme of proof of Theorem $1$.}

We fix a continuous map $f: S(X)\rightarrow \mathbb{R}^k$. The proof of theorem $1$ goes as follows, 
\begin{itemize}
\item Use a generalisation of the Borsuk-Ulam theorem giving rise to a finite convex partition of the sphere and a fiber of $f$ (i.e $f^{-1}(z)$ for some $z\in \mathbb{R}^k$) passing through the centers of all the pieces of the partition (the center of a convex set has to be defined).
\item Narrow the pieces of the partition (by increasing their numbers) such that almost all of them are Hausdorff close to a $k$-dimensional convex set. Pass to a limit infinite partition of the sphere by convex subsets of dimension less than or equal to $k$.
\item On each piece of the partition, there exists a probability measure, convexely derived from the conical measure. This brings the $n$-dimensional volume estimate of the waist down to a $k$-dimensional measure estimate on each convex set of the partition.
\\ This method is called the localization technique. But usually, the localization or the needle decomposition, brings the $n$-dimensional measure estimate down to a $1$-dimensional problem. The use of a multi-dimension localization technique first appears in \cite{grwst}.
\item On each piece of the partition, the Lemma $5.3$ gives an estimate of the measure of an $\varepsilon$-ball centered at a point where the measure of the convex set is mostly concentrated. By integrating this estimate over the space of pieces of the partition, we obtain the result of theorem $1$.
\\ There is some difficulties due to the $l$-dimensional convex sets of the infinite partition for all $l<k$. We prove that these "bad" convex sets does not affect the estimation of waist. Or better say, the measure of these convex sets in the space of pieces of the partition is equal to zero.
\end{itemize}

\section{Convexely derived measures on convex sets of $S(X)$}
The topics studied in this section follows the ideas used in \cite{ale} and \cite{gromil}.
For every subset $S \in S(X)$ we define the subset $co(S) \in B(X)$ as
$$co(S) :=\{\bigcup t S \vert  0 \leq t \leq 1\}.$$ 
Hence $co(S)$ is the cone centered at the origin of the ball over $S$.

\begin{de}[Convexely derived measure]
A convexely derived measure on $S(X)$ is a limit of a vaguely converging sequence of probability measures of the form $\mu_i=\frac{\mu|S_i}{\mu(S_i)}$, where $S_i$ are open convex sets.
\end{de}

Suppose we have a sequence of open convex sets $\{S_i\}$ of $S(X)$ which Hausdorff converges to a convex set $S' \in S(X)$ where we suppose that the dimension of $S'$ is equal to $k$ with $k<n$. It is clear that the sequence $\{co(S_i)\}$ Hausdorff converges to the set $co(S')$ where $\mathrm{dim}\, co(S')=k+1$. We define a probability measure $\mu'$ on $co(S')$ as follows.
For every $i \in \mathbb{N}$, we define the measure $\mu'_i=\frac{m_{n+1}|S_i}{m_{n+1}(S_i)}$. A subsequence of this sequence of measures vaguely converges to a probability measure $\mu$ on $co(S')$. We call this measure a convexely derived measure. We recall that the support of the measure $\mu$ is automatically equal to $co(S')$ as the sequence converges to this set. In \cite{ale} the author shows that the measure $\mu$ is $(n+1-(k+1))$-concave so by Borell's Theorem, $\mu$ admits a density function $f$ with respect to the $(k+1)$-dimensional Lebesgue measure defined on $A$. The function $f$ is $(n-k)$-concave. Hence
$$\mu= f dm_{k+1}$$ 
where $m_{k+1}$ is the $(k+1)$-dimensional Lebesgue measure. Morever we have:
\begin{lem}
 The measure $\mu$ is $(n+1)$-homogeneous and the function $f$ is $(n-k)$-homogeneous.
\end{lem}
This means $\mu(tA)=t^{n+1}\mu(A)$  for $0 \leq t \leq 1$
and $f(tx)=t^{n-k}f(x)$ for all $x \in co(S')$.

\emph{Proof of the Lemma}

The measure $\mu$ is convexely derived from the normalized $(n+1)$-dimensional Lebesgue measure. As the $(n+1)$-dimensional Lebesgue 
measure is $(n+1)$-homogeneous then $\mu$ is $(n+1)$-homogeneous.
From the equality $\mu=f dm_{k+1}$, and the fact that $\mu$ is $(n+1)$-homogeneous and $m_{k+1}$ is $(k+1)$-homogeneous, then clearly $f$
is $(n-k)$-homogeneous and the proof of the Lemma follows.

The convexely derived measure $\mu'$ defined on $co(S')$ defines a probability measure $\mu$ on $S'$ convexely derived from the conical measure of $S(X)$ and obtained from the sequence $\{S_i\}$, where for every $X\subset S'$ we have 
$$\mu(X)=\mu'(co(X)).$$
 And on the other hand, there exists another probability measure defined on $S'$ which is the canonical $k$-dimensional conical measure conically induced by $m_{k+1}$, we denote this measure by $\nu$. For every Borel subset $U$ of $S'$
$$\nu(U)=\frac{m_{k+1}(co(U))}{m_{k+1}(co(S'))}.$$
$S'$ is a subset of the unit sphere of $\mathbb{R}^{k+1}$ equipped with a norm satisfying the same modulus of convexity.

Then we have
$$\mu(U)=\mu'(co(U))=\int_{co(U)}f dm_{k+1}=\int_{U}f d\nu.$$
Hence in conclusion we have 
$$d\mu=f d\nu$$
where we take the restriction of $f$ on the set $U$.

The function $f$ is $(n-k)$-concave on $co(A)$ but the restriction of this function on the spherical part of the border of $co(A)$ is not anymore $(n-k)$-concave.

However the restriction function still has nice concavity properties as we will explain now.

\begin{de}
An {\em arc} $\sigma \subset S(X)$ is subarc of the intersection of a $2$-plane passing through the origin of the ball with $S(X)$.
\end{de}

We know that $\forall x,y \in S_{\pi}$,
$$f^{1/(n-k)}(\frac{x+y}{2})\geq \frac{f^{1/(n-k)}(x)+f^{1/(n-k)}(y)}{2}.$$

But the point $\frac{x+y}{2}$ is no more on $S(X)$, so we set 
$z=\frac{x+y}{2}/\Vert\frac{x+y}{2}\Vert \in S(X)$.

By the definition of the modulus of convexity we have
\begin{equation} \label{mod}
\Vert\frac{x+y}{2}\Vert\leq 1-\delta(\Vert x-y \Vert)
\end{equation}
So we can conclude the following Lemma.

\begin{lem} \label{conc}
Let $f$ denote the density of a convexely derived measure on $S(X)$.
Let $x$, $y \in S_{\pi}$, let $z=\frac{x+y}{2}/\Vert\frac{x+y}{2}\Vert \in S_{\pi}$. Then
$$\frac{f^{1/(n-k)}(x)+f^{1/(n-k)}(y)}{2}\leq(1-\delta(\Vert x-y \Vert)f^{1/(n-k)}(z).$$
\end{lem}

\emph{Proof of the Lemma}

As $\frac{x+y}{2}=\Vert\frac{x+y}{2}\Vert z$ and as the function $f$ is $(n-k)$-homogeneous 
$$f^{1/(n-k)}(\frac{x+y}{2})=\Vert\frac{x+y}{2}\Vert f^{1/(n-k)}(z)$$
and by equation \ref{mod} the proof of the Lemma follows.

\begin{de}
Let $f$ be a function defined on an arc of $S(X)$. Say $f$ is {\em weakly $(n-k)$-concave} if $\forall x$, $y \in \sigma$, $z=\frac{x+y}{2}/\Vert\frac{x+y}{2}\Vert$,
$$\frac{f^{1/(n-k)}(x)+f^{1/(n-k)}(y)}{2}\leq(1-\delta(\Vert x-y \Vert)f^{1/(n-k)}(z).$$ 
\end{de}

\begin{lem} \label{maxim}
A nonzero weakly $(n-k)$-concave function defined on an arc of $S(X)$ has at most one maximum point and has no local minima.
\end{lem}

\emph{Proof of the Lemma} 

If there were two distinct maxima $x$ and $y$ and we would get $f^{1/(n-k)}(x)\leq(1-\delta(\Vert x-y \Vert)f^{1/(n-k)}(x)$, contradiction. Suppose $f$ has a local minimum at point $m$. Take nearby points $x'$ and $y'$ such that $m=\frac{x'+y'}{2}$. Then $x=\frac{x'}{\Vert x'\Vert}$ and $y=\frac{y'}{\Vert y'\Vert}$ belong to the arc, and $m=\frac{x+y}{2}/\Vert\frac{x+y}{2}\Vert=m$. This leads again to a contradiction. The proof of the Lemma follows.

Let $f$ be the density of a convexely derived measure on supported on a $k$-dimensional convex subset $S$ of $S(X)$. By Lemma \ref{maxim} we can conclude that there exists at most one point $z \in S$ at which $f$ achieves its maximum. Indeed suppose the $f$ achieves its maximum in at least two points $x_1$ and $x_2$. Since there exists an arc passing through $x_1$ and $x_2$ and contained in $S$, this would contradict Lemma \ref{maxim}.

Let $z$ be the point of $S$ where $f$ achieves its maximum. We want to give a (uniform) lower bound for $\mu(B(z,\varepsilon))$ where $B(z,\varepsilon)$ is the $k$-dimensional ball in $S$ of norm-radius $\varepsilon$,
$$B(z,\varepsilon):=\{x\in S_{\pi}\vert \quad \Vert x-z \Vert \leq \varepsilon\}.$$
Therefore, from now on, the mm-space we are working on is $(S,\mu,\Vert \quad \Vert)$.

We define two subsets on $S$:
$A:=B(z,\varepsilon)$,
$B:=S\smallsetminus B(z,2\varepsilon)=B(z,2\varepsilon)^{c}$
and we are interested in estimating the ratio
$$
\frac{\mu_\pi(B)}{\mu_\pi(A)}.
$$

We need the following lemma.

\begin{lem} \label{comp}
Let $f$ be the density of a convexely derived measure supported on a $k$-dimensional convex subset $S$ of $S(X)$. Assume $f$ achieves its maximum at $z$. Let  $x \in B(z,2\varepsilon)^{c}=S\smallsetminus B(z,2\varepsilon)$ and consider the arc $\sigma=[z,x]$ in $S(X)$. Then
$$
f(x)\leq (1-2\delta(\varepsilon))^{n-k}Min\underset{\sigma\cap B(z,\varepsilon)}{f}.
$$
\end{lem}

\emph{Proof of the Lemma}

(Compare \cite{ale}) Pick $y \in [x,z]\cap B(z,\varepsilon)$. By weak concavity, we know that $f$ is monotone nondecreasing along $[x,z]$, so 
$$f(x)\leq f(y)\leq f(z).$$
So the maximum of $f$ on the subarc $[x,y]\sigma$ is achieved at $y$. By  Lemma \ref{conc},
$$\frac{f^{1/(n-k)}(x)+f^{\frac{1}{n-k}}(y)}{2}\leq(1-\delta(||x-y||)\underset{w\in [x,y]}{Max}f^{1/(n-k)}(w),$$
which implies
$$f(x)\leq(1-2(\delta(\Vert x-y\Vert))^{n-k} f(y).$$

By the triangle inequality, $\Vert x-y \Vert \geq \varepsilon$ and we remember that the modulus of convexity is nondecreasing, so 
$$
\delta(\Vert x-y \Vert)\geq \delta(\varepsilon).
$$
Hence
$$
(1-2\delta(\Vert x-y \Vert))^{n-k}\leq (1-2\delta(\varepsilon))^{n-k}.
$$
And at last we have
$$
f(x)\leq (1-2\delta(\Vert x-y \Vert))^{n-k}f(y)\leq (1-2\delta(\varepsilon))^{n-k}f(y).
$$
And the proof of the Lemma follows.

We are ready now to integrate both sides of the inequality of Lemma \ref{comp} and give an upper bound for $\frac{\mu(B)}{\mu(A)}.$

\begin{lem}
Let $\varepsilon>0$ be given. Let $S\subset S(X)$ be a $k$-dimensional convex set. Let a convexely derived measure $\mu$ be defined on $S$. Let $z$ be the maximum point for the density function of the measure $\mu$. Let $A:=B(z,\varepsilon)$, $B:=S\smallsetminus B(z,2\varepsilon)$. Then
\begin{eqnarray*}
\frac{\mu(B)}{\mu(A)}&\leq&(1-2\delta(\varepsilon))^{n-k}(k+1)^{k+1}\frac{F(k,\varepsilon)}{G(k,\varepsilon)}.
\end{eqnarray*}
$$F(k,\varepsilon)= \int_{\psi_2(\varepsilon)}^{\frac{\pi}{2}}\sin(x)^{k-1}\,dx.$$ 
and 
$$G(k,\varepsilon)= \int_{0}^{\psi_1(\varepsilon)}\sin(x)^{k-1}\,dx.$$
And where 
$$\psi_1(\varepsilon)=2\arcsin(\frac{\varepsilon}{4\sqrt{k+1}})$$
and
$$\psi_2(\varepsilon)=2\arcsin(\frac{\varepsilon}{2\sqrt{k+1}})$$
\end{lem}

\emph{Proof of the Lemma}

Let $\sigma$ be an arc of $S(X)$ emanating from $z$. Denote by
$$
m=Min\underset{\sigma\cap B(z,\varepsilon)}{f}.
$$
Then
$$
x \in \sigma\cap B(z,2\varepsilon)^{c} \Rightarrow f(x)\leq (1-2\delta(\varepsilon))^{n-k}m ,
$$
and
$$
y \in \sigma\cap B(z,\varepsilon) \Rightarrow f(y)\geq m.
$$

Assume first that the norm $\Vert\cdot\Vert$ is Euclidean. We need to convert Euclidean distances into Riemannian distances along the unit sphere, i.e. angles. If $x$ and $y$ are unit vectors making an angle $\phi$, then 
$|x-y|=2\sin(\phi/2)$. Therefore $|x-y|=\epsilon$ corresponds to an angle $\phi_1$ and $|x-y|=2\epsilon$ corresponds to an angle $\phi_2$. Therefore, for a fixed $\theta$, $t\leq\phi_1 \Rightarrow f(t,\theta)\geq m(\theta)$ and 
$t\geq \phi_2 \Rightarrow f(t,\theta)\leq (1-2\delta(\varepsilon))^{n-k}m(\theta) $. 
Using polar coordinates $(t,\theta)$ on the unit sphere, we compute
\begin{eqnarray*}
\frac{\mu(B)}{\mu(A)}&\leq&\frac{\int_{\phi_2}^{\pi} \int_{\mathbb{S}^{k-1}}f(t,\theta)\,\sin(t)^{k-1}\,dt\,d\theta}{\int_{0}^{\phi_1} \int_{\mathbb{S}^{k-1}}f(t,\theta)\,\sin(t)^{k-1}\,dt\,d\theta}\\
&\leq&\max_{\theta\in\mathbb{S}^{k-1}}\frac{\int_{\phi_2}^{\pi} f(t,\theta)\,\sin(t)^{k-1}\,dt}{\int_{0}^{\phi_1}f(t,\theta)\,\sin(t)^{k-1}\,dt}.
\end{eqnarray*}
For each $\theta$,
\begin{eqnarray*}
\frac{\int_{\phi_2}^{\frac{\pi}{2}} f(t,\theta)\,\sin(t)^{k-1}\,dt}{\int_{0}^{\phi_1}f(t,\theta)\,\sin(t)^{k-1}\,dt}&\leq&\frac{\int_{\phi_2}^{\pi}(1-2\delta(\varepsilon))^{n-k}m(\theta)\sin(t)^{k-1}\,dt}{\int_{0}^{\phi_1}m(\theta)\sin(t)^{k-1}\,dt}\\
&=&\frac{\int_{\phi_2}^{\pi}\sin(t)^{k-1}\,dt}{\int_{0}^{\phi_1}\sin(t)^{k-1}\,dt}(1-2\delta(\varepsilon))^{n-k}.
\end{eqnarray*}

To handle general norms, we use the fact that the Banach-Mazur distance between any $k+1$-dimensional normed space and Euclidean space is at most $\sqrt{k+1}$. On the affine extension of $co(S)$ there exists a Euclidean structure $\vert \cdot \vert$ such that for every $x \in Aff(co(S))$ we have
$$
\frac{1}{\sqrt{k+1}}\vert x \vert \leq \Vert x \Vert \leq \vert x \vert.
$$
Or equivalently we have
$$B \subset K \subset \sqrt{k+1}B,$$
where $B$ is the Euclidean ball of dimension $k+1$ and $K$ is the uniformly convex ball defined by $S(X)$.

We denote by $pr$ the radial projection of the uniformly convex sphere $\partial K$ to the Euclidean sphere $\partial B$. Recall that $\nu$ is the conical measure on $\partial K$ and we denote by $dv_{k}$ the conical measure on $\partial B$, i.e. the Riemannian probability measure. Then the density $\displaystyle h=\frac{pr_{*}d\nu}{dv_{k}}$ satisfies
\begin{eqnarray*}
\frac{1}{\sqrt{k+1}^{k+1}}\leq h \leq \sqrt{k+1}^{k+1} .
\end{eqnarray*}

Let $x$, $y\in\partial K$, $x'=pr(x)$, $y'=pr(y)$. Since radial projection to the sphere decreases Euclidean distance outside the Euclidean ball,
\begin{eqnarray*}
|x'-y'|\leq |x-y|\leq \sqrt{k+1}\Vert x-y \Vert.
\end{eqnarray*}
For a general norm, radial projection to the unit sphere is $2$-Lipschitz. Indeed, let $x''$, $y''$ be points such that $1\leq\Vert x''\Vert\leq\Vert y''\Vert$. Rescaling both by $\Vert x''\Vert$ decreases $\Vert x''-y''\Vert$, so we can assume that $\Vert x''\Vert=1$. Then $\Vert y''\Vert\leq 1+\Vert x''- y''\Vert$ and
\begin{eqnarray*}
\Vert x''-\frac{y''}{\Vert y''\Vert}\Vert
&=&\Vert \frac{x''}{\Vert y''\Vert}-\frac{y''}{\Vert y''\Vert}+x(1-\frac{1}{\Vert y''\Vert})\Vert\\
&\leq&\Vert x''-y''\Vert+\Vert y''\Vert-1\leq 2\Vert x''-y''\Vert.
\end{eqnarray*}
If $x''=\sqrt{k+1}x'$ and $y''=\sqrt{k+1}y'$, then
\begin{eqnarray*}
\Vert x-y \Vert\leq 2\Vert x''-y'' \Vert=2\sqrt{k+1}\Vert x'-y' \Vert\leq 2\sqrt{k+1}|x'-y'|.
\end{eqnarray*}

We radially project the set $S$ to a set $S'$ on the sphere. $S'$ is $k$ dimensional and is a convex set as radial projection preserves convexity. We denote the projection of the point $z$ on the sphere by $z'=pr(z)$. In polar coordinates $(t,\theta)$ centered at $z'$, fix $\theta$. Let $\psi_1 (\theta)$ (resp. $\psi_2 (\theta)$) denote the angle $t$ such that $y=pr^{-1}(t,\theta)\in\partial K$ satisfies $\Vert y-z\Vert=\varepsilon$ (resp. $=2\varepsilon$). The above distance estimates yield
\begin{eqnarray*}
2\sin\frac{\psi_1 (\theta)}{2}\geq \frac{\varepsilon}{2\sqrt{k+1}}
\end{eqnarray*}
and
\begin{eqnarray*}
2\sin\frac{\psi_2 (\theta)}{2}\geq \frac{\varepsilon}{\sqrt{k+1}}.
\end{eqnarray*}
Then
\begin{eqnarray*}
\frac{\mu(B)}{\mu(A)}&\leq&\frac{\int_{\mathbb{S}^{k-1}}\int_{\psi_2 (\theta)}^{\pi} h(t,\theta)f(t,\theta)\,\sin(t)^{k-1}\,dt\,d\theta}{\int_{\mathbb{S}^{k-1}}\int_{0}^{\psi_1 (\theta)} h(t,\theta)f(t,\theta)\,\sin(t)^{k-1}\,dt\,d\theta}\\
&\leq&\max_{\theta\in\mathbb{S}^{k-1}}\frac{\int_{\psi_2 (\theta)}^{\pi} h(t,\theta)f(t,\theta)\,\sin(t)^{k-1}\,dt}{\int_{0}^{\psi_1 (\theta)}h(t,\theta)f(t,\theta)\,\sin(t)^{k-1}\,dt}.
\end{eqnarray*}
For each $\theta$,
\begin{eqnarray*}
\frac{\int_{\psi_2 (\theta)}^{\pi} h(t,\theta)f(t,\theta)\,\sin(t)^{k-1}\,dt}{\int_{0}^{\psi_1 (\theta)}h(t,\theta)f(t,\theta)\,\sin(t)^{k-1}\,dt}&\leq&\frac{\int_{\psi_2}^{\pi}(1-2\delta(\varepsilon))^{n-k}m(\theta)h(t,\theta)\sin(t)^{k-1}\,dt}{\int_{0}^{\psi_1}m(\theta)h(t,\theta)\sin(t)^{k-1}\,dt}\\
&=&(1-2\delta(\varepsilon))^{n-k}\frac{\int_{\psi_2}^{\pi}h(t,\theta)\sin(t)^{k-1}\,dt}{\int_{0}^{\psi_1}h(t,\theta)\sin(t)^{k-1}\,dt}\\
&\leq&(1-2\delta(\varepsilon))^{n-k}(k+1)^{k+1}\frac{\int_{\psi_2}^{\pi}\sin(t)^{k-1}\,dt}{\int_{0}^{\psi_1}\sin(t)^{k-1}\,dt}.
\end{eqnarray*}
Replacing $\psi_1$ and $\psi_2$ with the above lower bounds yields
\begin{eqnarray*}
\frac{\mu(B)}{\mu(A)}&\leq&(1-2\delta(\varepsilon))^{n-k}(k+1)^{k+1}\frac{\int_{\psi_2}^{\pi}\sin(t)^{k-1}\,dt}{\int_{0}^{\psi_1}\sin(t)^{k-1}\,dt}\\
                     &\leq&(1-2\delta(\varepsilon))^{n-k}(k+1)^{k+1}\frac{F(k,\varepsilon)}{G(k,\varepsilon)}.
\end{eqnarray*}

And the proof of the Lemma follows.
\begin{lem}
Let $S$ be a convex set of dimension $k$ in $S(x)$. Let a convexely derived measure $\mu$ be defined on $S$. Let $z$ be the maximum point of the density of the measure $\mu$. For every $\varepsilon>0$ we have the following estimation
$$\mu(B(z,\varepsilon)\geq \frac{1}{1+(1-2\delta(\frac{\varepsilon}{2}))^{n-k}(k+1)^{k+1}\frac{F(k,\frac{\varepsilon}{2})}{G(k,\frac{\varepsilon}{2})}}$$

Where the functions $F$ and $G$ are as defined before.
\end{lem}
\emph{Proof of the Lemma}
We use the result of the previous Lemma which tells 
\begin{eqnarray*}
\frac{\mu(B)}{\mu(A)}&\leq&(1-2\delta(\varepsilon))^{n-k}(k+1)^{k+1}\frac{F(k,\varepsilon)}{G(k,\varepsilon)}.
\end{eqnarray*}
We remind that $\mu$ is a probability measure and we have
$$\frac{\mu(B(z,2\varepsilon))}{\mu(B(z,2\varepsilon))^{c}}\geq\frac{\mu(B(z,\varepsilon))}{\mu(B(z,2\varepsilon))^{c}}\geq \frac{1}{(1-2\delta(\varepsilon))^{n-k}(k+1)^{k+1}\frac{F(k,\varepsilon)}{G(k,\varepsilon)}}.$$
Hence
$$
\mu(B(z,2\varepsilon))=\frac{\mu(B(z,2\varepsilon))}{\mu(B(z,2\varepsilon))+\mu(B(z,2\varepsilon))^{c}}\geq \frac{1}{1+(1-2\delta(\varepsilon))^{n-k}(k+1)^{k+1}\frac{F(k,\varepsilon)}{G(k,\varepsilon)}}
$$
And the proof of the Lemma follows.

\section{Proof of Theorem 1 following Gromov}
In this section we follow the ideas used in \cite{grwst} and \cite{memwst}. Let $f:S(X) \to \mathbb{R}^k$ be as theorem $1$. We want to partition the sphere $S(X)$ by at most $k$-dimensional convex sets. The continuous map $f$ defines a continuous map $Pr_{r}(f)$ on the sphere $\mathbb{S}^n$ which is the radial projection of $f$ on $\mathbb{S}^n$. We use the following theorem proved announced by Gromov in \cite{grwst}. The author remarks that the following Theorem is not entirely proved in \cite{grwst} and unfortunately we are not able to give a proof for this Theorem. however, if we believe Gromov, then the proof of our Theorem $1$ becomes much easier. In the other hand, we will give another method, which will be independent of the following Theorem to finalize the results of this paper.
\begin{theo}[Gromov] \label{infi}
Let $f : \mathbb{S}^n \rightarrow \mathbb{R}^k$ be a continuous map. There exists an infinite partition of the sphere by at most $k$-dimensional convex sets, denoted by $\Pi_{\infty}$ and a point $z \in \mathbb{R}^k$ such that for every $S \in \Pi_{\infty}$, $f^{-1}(z)$ passes through the maximum point of the density of the convexely derived measure defined on $S$.
\end{theo}
We inform the reader that the previous Theorem holds for every continuous map $f$.
\begin{cor}
Let $f: S(X) \rightarrow \mathbb{R}^k$ be as theorem $1$. There exist an infinite partition of $S(X)$ by at most $k$-dimensional convex sets, denoted by $\Pi_{\infty}$ and a point $z \in \mathbb{R}^k$ such that for every $S \in \Pi_{\infty}$, $f^{-1}(z)$ passes through the maximum point of the density of the (unique) convexely derived measure defined on $S$.
\end{cor}
\emph{Proof of the Corollary}
 we apply the previous theorem for the continuous map $Pr(f)$, we know that there exists an infinite partition of the sphere, $\Pi_{\infty}$, by at most $k$-dimensional convex sets. By radialy projecting each piece of the partition on $S(X)$ we obtain an infinite partition of $S(X)$ by at most $k$-dimensional convex sets. Let $S \subset \mathbb{S}^n$ and $S \in \Pi_{\infty}$ and let $S'=pr(S)$. Denote by $z$ (resp $z'$) the maximum point of the density of the convexely derived measure defined on $S$ (resp $S'$). It remains to prove that $z'=pr(z)$. Indeed as we are taking radial projection, the density of the convexely derived mesure on each $S'$ is just the radial projection of the density of the measure defined on $S$. We remind that the radial projection of the normalized Riemannian measure of $S$ is the conical measure defined on $S'$ up to a constant, but this is irelevant for our purpose.
We are ready to give a proof of the Theorem $1$.
\subsection{Proof of Theorem 1 following Theorem \ref{infi}}
We apply the previous Corollary. There exists an infinite partition of $S(X)$ by at most $k$-dimensional convex sets and a fiber $f^{-1}(z)$ passing through all the maximum points of the densities of the convexely derived mesure defined on all pieces of the partition.
where $x_{\pi}$ is the maximum point of the density of the (unique) convexely derived measure $\mu_{\pi}$ defined on $S_{\pi}$.
Hence on every $S_{\pi}$ we have
$$\mu_{\pi}((f^{-1}(z)+\varepsilon)\cap S_{\pi})\geq \mu_{\pi}(B(x_{\pi},\varepsilon))\geq w(\varepsilon).$$
And at the end
$$\mu(f^{-1}(z)+\varepsilon)=\int_{\Pi_{\infty}}\mu_{\pi}((f^{-1}(z)+\varepsilon)\cap S_{\pi})d\pi$$
$$=\int_{dim S_{\pi}=k}\mu_{\pi}((f^{-1}(z)+\varepsilon)\cap S_{\pi})d\pi+\int_{dim S_{\pi}<k}\mu_{\pi}((f^{-1}(z)+\varepsilon)\cap S_{\pi})d\pi$$
The measure of the measurable partition is equal to one. In \cite{memwst} we prove that the measure of the set of pieces of partition which has dimension $<k$ on the sphere is equal to zero, radialy projecting this on $S(X)$ implies that the measure of the set of pieces of partition of $S(X)$ which has dimension $<k$ is also equal to zero, hence we have
$$
\mu(f^{-1}(z)+\varepsilon)\geq w(\varepsilon).
$$
Hence the proof of the theorem follows.
\section{Alternative proof of Theorem 1}
\label{infinite}
This section will be long and very technical. As the author is unable to prove Theorem \ref{infi}, he found, by the enormous help of Pierre Pansu, the following arguments replacing theorem \ref{infi}. We also remark that the obstruction for having theorem \ref{infi} is due to non-existence of a sharp Brunn-Minkowski type inequality on the round sphere. We begin by giving the following useful
\begin{de}
Let $S$ be an open convex subset of $S(X)$, $S$ is called an $(k,\varepsilon)$-pancake if there exists a convex set $S_\pi$ of dimension $k$ such that every point of $S$ is at distance at most $\varepsilon$ from $S_\pi$.
\end{de}
We remark again that the distance on $S(X)$ is the restriction of the norm being defined on $\mathbb{R}^{n+1}$ on $S(X)$.

The two following theorems are strong generalizations of the classical Borsuk-Ulam theorem in algebraic topology and the construction of finite and infinite partitions of $S(X)$ is provided by them.

\begin{theo}[Gromov-Borsuk-Ulam, finite case]
\label{3}
Let $f: \mathbb{S}^{n} \rightarrow \mathbb{R}^k$ ($k\leq n$) be a continuous map from the $n$-sphere to Euclidean space of dimension $k$. For every $i \in  \mathbb{N}$, there exists a partition of the sphere $\mathbb{S}^n$ into $2^i$ open convex sets $\{S_i\}$ of equal volumes ($=Vol(S^n)/2^i$) and such that all the center points $c_{.}(S_i)$ of the elements of partition have the same image in $\mathbb{R}^k$.
\end{theo} 

\begin{theo}[Gromov-Borsuk-Ulam, almost infinite case]
\label{4}
Let $f:\mathbb{S}^n \rightarrow \mathbb{R}^k$ be a continuous map. For all $\varepsilon>0$, there exists an integer $i_0$ such that for all $i \geq i_0$ there exists a finite partition of $\mathbb{S}^n$ into $2^i$ open convex subsets such that :
\begin{enumerate}
\item Every convex subset of the partition is a $(k,\varepsilon)$-pancake.
\item The centers of all convex subsets of the partition have the same image in $\mathbb{R}^k$.
\item All convex subsets of the partition have the same volume.
\end{enumerate}
\end{theo}
The proof of theorem \ref{3} is long and uses algebraic topology arguments. We won't give the proof of these theorems here and refer the reader to \cite{memwst}. 

We need Theorems \ref{3} and \ref{4} on $S(X)$, but we can not proceed directly, we again pass via the round sphere and by radially projecting the results of these two theorems on $S(X)$, we obtain the desired partitions on $S(X)$.

\subsection{Approximation of General Norms By Smooth Norms}
For technical reason imposed by Lemma \ref{radon}, we need to approximate general norms by smooth norms. Indeed as we will see in the next subsection, we can not allow the convexely derived measures charging any mass for the boundary of balls. In this subsection we show by approximation that we can in fact exclude this technical problem.
\begin{lem}
\label{courbure}
Let $X$ denote a finite dimensional space equipped with a $C^2$-smooth norm. Let $S(X)$ denote its unit sphere. Fix an auxiliary Euclidean structure. There exists $K$ such that for every $2$-plane $\Pi$ passing through the origin, $S(X)\cap P$ is a disjoint union of curves whose curvatures $\kappa$ satisfy $|\kappa|\leq K$ at all points.
\end{lem}

Proof.

Since the norm is homogeneous of degree $1$, its derivative along a line passing through the origin does not vanish. It follows that at every point $x\in S(X)$, the restriction of the differential to $P$ does not vanish identically, i.e. $P$ is transverse to the tangent hyperplane $T_x S(X)$. This shows that $S(X)\cap P$ is a $C^2$-smooth $1$-dimensional submanifold, i.e. a finite disjoint union of curves. Furthermore, the curvature $\kappa(x,P)$ of $S(X)\cap P$ at $x$ is a continuous function of $(x,P)\in I=\{(x,P)\,|\,x\in\partial B(0,1),\,x\in P\}$. Since $I$ is compact, $\kappa$ is bounded.

\begin{nott}
\label{hessian}
The {\em Hessian} of a $C^2$-smooth function $f:\mathbb{R}^{d}\to\mathbb{R}$ at $x$ is the quadratic form 
\begin{eqnarray*}
Hess_{x}(v)=\frac{\partial^{2}}{\partial t^{2}}f(x+tv)_{|t=0}.
\end{eqnarray*}
Say a $C^2$-smooth norm on a finite dimensional vectorspace is {\em strongly convex} if at every nonzero point, the Hessian of $x\mapsto \parallel x\parallel^2$ is positive definite.
\end{nott}

\begin{prop}
\label{arcs}
Let $X$ denote a finite dimensional space equipped with a $C^2$-smooth strongly convex norm. Let $S(X)$ denote its unit sphere. There exists $r_0 >0$ such that, for every $r<r_0$, for every $2$-plane $P$ passing through the origin, for every $x\in S(X)$, $S(X)\cap P\cap\partial B(x,r)$ is a finite set.
\end{prop}

Proof.

The map $x\mapsto Hess_{x}\parallel \cdot\parallel^2$ is homogeneous of degree $0$. Fix an auxiliary Euclidean inner product on $X$. By compactness of the unit sphere, there exists a positive constant $c$ such that for all $x\not=0$ and all $v$,
\begin{eqnarray}
(Hess_{x}\parallel \cdot\parallel^2) (v,v)\geq c\,v\cdot v. \label{hess}
\end{eqnarray}
Also, the differential $x\mapsto D_{x}\parallel \cdot\parallel^2$ is homogeneous of degree $1$. Therefore there exists a positive constant $C$ such that for all $x\not=0$ and all $v$,
\begin{eqnarray}
|(D_{x}\parallel \cdot\parallel^2)(v)|\leq C\,\parallel x\parallel\sqrt{v\cdot v}. \label{diff}
\end{eqnarray}
Fix $x\in X$. Let $P$ be a $2$-plane. Let $f$ denote the restriction of $z\mapsto\parallel z-x\parallel^2$ to $P$. It satisfies the previous two inequalities. Let $s\mapsto \gamma(s)$ be a $C^2$-smooth curve in $P$ parametrized by arclength, $z=\gamma(0)$, $\tau=\gamma'(0)$. Then
\begin{eqnarray*}
\gamma(s)=z+s\tau+\frac{s^2}{2}\gamma''(0) +o(s^2).
\end{eqnarray*}
Since, for all small $v$,
\begin{eqnarray*}
f(z+v)=f(z)+D_{z}f(v)+\frac{1}{2}Hess_{x}f(v,v)+o(v\cdot v),
\end{eqnarray*}
\begin{eqnarray*}
f(\gamma(s))=f(z)+D_{z}f(s\tau+\frac{s^2}{2}\gamma''(0))+\frac{1}{2}Hess_{z}f(\tau,\tau)+o(s^2).
\end{eqnarray*}
Now assume that $f(\gamma(s_j))=f(z)$ for a sequence $s_j$ that tends to $0$. Then, comparing asymptotic expansions gives
\begin{eqnarray*}
D_{z}f(\tau)=0,\quad D_{z}f(\gamma''(0))+Hess_{z}f(\tau,\tau).
\end{eqnarray*}
Since $\tau\cdot\tau=1$, inequalities (\ref{hess}) and (\ref{diff}) give
\begin{eqnarray*}
c\leq -D_{z}f(\gamma''(0))\leq C\,\parallel z-x\parallel\sqrt{\gamma''(0)\cdot \gamma''(0)}.
\end{eqnarray*}
This shows that the curvature $\kappa$ of the plane curve at $\gamma$ at $z$ satisfies
\begin{eqnarray*}
\kappa(z)\geq\frac{c}{C\parallel z-x\parallel}.
\end{eqnarray*}
Therefore, if $z$ is an accumulation point of $\gamma\cap P\cap\partial B(x,r)$, the curvature of $\gamma$ ay $z$ is $\geq \frac{c}{Cr}$. With Lemma \ref{courbure}, we conclude that if $r<r_0 :=c/CK$, for all $P$, $S(X)\cap P\cap\partial B(x,r)$ has only isolated points, thus is finite.

\begin{lem}
\label{approx}
Let $X_1$ be a finite dimensional normed space. Let $S(X_1)$ denote its unit sphere. For every $\lambda>1$, there exists a $C^2$-smooth strongly convex norm on $X_1$, with unit sphere $S(X_2)$, such that the radial projection $S(X_1)\to S(X_2)$ is $\lambda$-biLipschitz.
\end{lem}

Proof.

Fix an auxiliary Euclidean inner product on $X_1$. Fix a smooth compactly supported nonnegative function $\psi:X\to\mathbb{R}_+$ such that $\int\psi=1$. The convolution
\begin{eqnarray*}
f(x)=\int_{X_1}\parallel y\parallel_1 \psi(x-y)\,dy=\int_{X_1}\parallel x-y\parallel_1 \psi(y)\,dy
\end{eqnarray*} 
is smooth and convex. For all $x\in X_1$,
\begin{eqnarray*}
|f(x)-\parallel x\parallel_1 |\leq\int_{X_1}\parallel y\parallel_1 \psi(y)\,dy
\end{eqnarray*}
is uniformly bounded. Therefore, when one restricts $f$ to a large Euclidean sphere and extends it to become positively homogeneous of degree $1$, one gets a smooth norm $\parallel \cdot\parallel'$ uniformly close to $\parallel \cdot\parallel_1$. By convexity, the Hessian of $\parallel \cdot\parallel'^2$ is nonnegative. For $\delta>0$, let
\begin{eqnarray*}
\parallel v\parallel_\delta =\sqrt{\parallel v\parallel'^{2}+\delta\,v\cdot v}.
\end{eqnarray*}
This is a smooth norm, and $Hess(\parallel v\parallel_{\delta}^{2})\geq\delta\,v\cdot v$ is positive definite. For $\delta$ small enough, this norm is close to $\parallel \cdot\parallel_1$, therefore radial projection between unit spheres is $\lambda$-biLipschitz.

\bigskip

Lemma \ref{approx} allows to reduce the proof of Theorem \ref{usphere} to the special case of $C^2$-smooth strongly convex norms, for which we know, from Proposition \ref{arcs}, that convexely derived measures do not give any mass to small enough spheres. Until the end of section \ref{infinite}, we suppose the norm of class $C^2$ and strongly convex.

\subsection{Infinite Partitions}
The proof follows \cite{memwst}, where the case of the round sphere $\mathbb{S}^n$ is treated. But we need these results for the unit spheres of uniformly convex normed spaces. This merely requires a few minor changes, but we include complete proofs for completeness sake.

\begin{de}[space of convexely derived measures]

Let $\mathcal{MC}^{n}$ de\-no\-te the set of probability measures on $S(X)$ of the form $\mu_S =\mu_{|S}/\mu(S)$ where $S\subset S(X)$ is open and convex and where $\mu$ is the conical probability measure defined on $S(X)$. The space $\mathcal{MC}$ of {\em convexely derived probability measures} on $S(X)$ is the vague closure of $\mathcal{MC}^{n}$.
\end{de} 

It is a compact metrizable topological space.

\begin{lem}
\label{bishop}
For all open convex sets $S\subset\mathbb{S}^n$ and all $x\in S$,
\begin{eqnarray*}
\frac{vol(S\cap B(x,r))}{vol(S)}\geq \frac{vol(B(x,r))}{vol(\mathbb{S}^n)}.
\end{eqnarray*}
\end{lem}

Proof.

Apply Bishop-Gromov's inequality in Riemannian geometry. In this special case ($\mathbb{S}^n$ has constant curvature $1$), it states that the ratio
\begin{eqnarray*}
\frac{vol(S\cap B(x,r))}{vol(B(x,r))}
\end{eqnarray*}
is a nonincreasing function of $r$. It follows that
\begin{eqnarray*}
\frac{vol(S\cap B(x,r))}{vol(B(x,r))}\geq\frac{vol(S)}{vol(\mathbb{S}^n)}.
\end{eqnarray*}
\begin{cor}
\label{bishoppp}
For all open convex sets $S\subset S(X)$ and all $x\in S$,
\begin{eqnarray*}
\frac{\mu(S\cap B(x,r))}{\mu(S)}\geq \frac{vol(.,\phi(r))}{vol(\mathbb{S}^n)},
\end{eqnarray*}
where $vol(.,\phi(r))$ is the volume of a ball of radius $\phi(r)$ on $\mathbb{S}^n$ and where 
\begin{eqnarray*}
2\sin(\frac{\phi(r)}{2})= \frac{r}{2\sqrt{n+1}}.
\end{eqnarray*}
\end{cor}

Proof.

By radially projecting $S(X)$ to $\mathbb{S}^n$, the convex set $S$ maps to a convex set $S'$ on the round sphere. By our previous observations, the image of the ball $B(x,r)$ contains a spherical ball of radius $\phi(r)$ where
\begin{eqnarray*}
2\sin(\frac{\phi(r)}{2})=\frac{r}{2\sqrt{n+1}}.
\end{eqnarray*}
Hence
\begin{eqnarray*}
\frac{\mu(S\cap B(x,r))}{\mu(S)}&\geq& \frac{\mu'(S'\cap B(x',r'))}{\mu'(S')}\\
&\geq&\frac{\mu'(S'\cap B(x',\phi(r)))}{\mu'(S')}\\
&\geq&\frac{(n+1)^{n+1}vol(B(x',\phi(r))}{(n+1)^{n+1}vol(\mathbb{S}^n)}\\
&=&\frac{vol(B(x',\phi(r)))}{vol(\mathbb{S}^n)}.
\end{eqnarray*}

This inequality extends to all convexely derived measures, thanks to the following Lemma.

\begin{lem}
\label{radon}
{\em (See \cite{hirsh}).}
Let $\mu_i$ be a sequence of positive Radon measures on a locally compact  space $X$ which vaguely converges to a positive Radon measure $\mu$. Then for every relatively compact subset $A \subset X$ such that $\mu(\partial A)=0$,
$$\lim_{i \to \infty}\mu_i (A)=\mu(A).$$
\end{lem}

\begin{cor}
\label{bishopcdm}
For all measures $\nu\in\mathcal{MC}$ on $S(X)$, all $x\in \mathrm{support}(\nu)$ and small enough $r$
\begin{eqnarray*}
\nu(S\cap B(x,r))\geq const.\,r^n .
\end{eqnarray*}
\end{cor}

Proof.

Let $\nu=\lim \mu_{S_{j}}$. Up to extracting a subsequence, one can assume that $S_{j}$ Hausdorff converges to a compact convex set $S$. Then $\mathrm{support}(\nu)\subset S$. Indeed, if $x\notin S$, there exists $r>0$ such that $S\cap B(x,r)=\emptyset$. Let $f$ be a continuous function on $S(X)$, supported in $B(x,r/2)$. Then for $j$ large enough, $S_{j}\cap B(x,r/2)=\emptyset$, $\int f\,d\nu_{S_{j}}=0$, so $\int f\,d\nu=0$, showing that $x\notin\mathrm{support}(\nu)$.

If $\nu$ is a Dirac measure, then the inequality trivially holds. Otherwise, let $x\in\mathrm{support}(\nu)$. There exist $x_j\in \mathrm{support}(\mu_j)$ such that $x_j$ tend to $x$. Since $\nu$ gives no measure to boundaries of small metric balls (by Proposition \ref{arcs}, since we assume that the norm is $C^2$ and strongly convex), Lemma \ref{radon} applies, and the inequality of Corollary \ref{bishopcdm} passes to the limit.

\begin{lem}
\label{haus}
Let $Comp(S(X))$ denote the space of compact subsets of $S(X)$ equipped with Hausdorff distance. The map $\mathrm{support}:\mathcal{MC}\to Comp(S(X))$ which maps a measure to its support is continuous.
\end{lem}

Proof.

Let $\mu_j \in\mathcal{MC}$ converge to $\nu$. One can assume that $S_j =\mathrm{support}(\mu_j)$ converge to a compact set $S$. We saw in the proof of Corollary \ref{bishopcdm} that $\mathrm{support}(\nu)\subset S$.
To prove the opposite inclusion, let us define, for $r>0$ and $x\in S(X)$,
\begin{eqnarray*}
f_{r,x}(y)=\begin{cases}
1 & \text{ if }d(y,x)<\frac{r}{2}, \\
2- 2\frac{d(y,x)}{r} & \text{ if }\frac{r}{2}\leq d(y,x)<r, \\
0 & \text{otherwise}.
\end{cases}
\end{eqnarray*}
Where $d$ is the distance induced by the norm of $\mathbb{R}^{n+1}$. Let $x\in S$. Let $x_j \in S_j$ converge to $x$. According to Lemma \ref{bishopcdm}, if  $d(x_j ,x)<r/4$,
$$
\int f_{x,r}(y)\,d\mu_{j}(y)\geq const.r^{n},
$$
i.e. $\displaystyle \int f_{x,r}\,d\mu_{j}$ does not tend to $0$. It follows that $\displaystyle \int f_{x,r}\,d\nu>0$, and $x$ belongs to $\mathrm{support}(\nu)$. This shows that $\mathrm{support}$ is a continuous map on $\mathcal{MC}$.

\medskip

The support of a convexely derived probability measure is a closed convex set, it has a dimension.

\begin{nott}
$\mathcal{MC}^{k}$ denotes the set of convexely derived probability measures whose support has dimension $k$, $\mathcal{MC}^{\leq k}=\bigcup_{\ell=0}^{k}\mathcal{MC}^{k}$, $\mathcal{MC}^{+}=\mathcal{MC}\setminus \mathcal{MC}^{0}$. For $\rho>0$, $\mathcal{MC}_{\rho}$ denotes the set of convexely derived probability measures whose support has diameter $\geq\rho$.
\end{nott}

\begin{lem}
\label{maj}
As $r$ tends to $0$, $\nu(B(x,r))$ tends to $0$ uniformly on $\mathcal{MC}_{\rho}\times\ S(X)$.
\end{lem}

Proof.

We first prove the Lemma in $\mathbb{R}^n$, the spherical case follows by projectively mapping hemispheres of $S(X)$ to $\mathbb{R}^n$. We can assume that $\rho$ is very small as well. Let $\mu$ be a convexely derived measure supported by a $k$-dimensional convex set $S$, let $x\in\mathbb{R}^n$ and $B=S\cap B(x,r)$. Since $S$ has diameter at least $\rho$, there is a point $y$ at distance at least $\rho/2$ of $x$. Up to a translation, we can assume that $y$ is the origin of $\mathbb{R}^k$. Let $\phi$ be the density of $\mu$. Then $\phi^{1/(n-k)}$ is concave. Thus, for $x'\in B$ and $\lambda\in]0,1[$,
\begin{eqnarray*}
\phi(\lambda x)\geq \lambda^{n-k}\phi(x).
\end{eqnarray*}
Changing variables gives
\begin{eqnarray*}
\mu(\lambda B)&=&\int_{\lambda B}\phi(z)\,dz\\
&=&\lambda^{k}\int_{B}\phi(\lambda z)\,dz\\
&\geq&\lambda^{n}\int_{B}\phi(z)\,dz\\
&=&\lambda^{n}\mu(B).
\end{eqnarray*}
If $N$ is an integer such that $N\leq \rho/4r$, then one can choose $N$ values of $\lambda$ between $1/2$ and $1$ leading to disjoint subsets $\lambda B$ of $S$, and this yields
\begin{eqnarray*}
1=\mu(S)\geq N(\frac{1}{2})^{n}\mu(B),
\end{eqnarray*}
i.e.
\begin{eqnarray*}
\mu(B)\leq 2^{n}/N \simeq \mathrm{const.}\,r/\rho.
\end{eqnarray*}
Now, let $S \subset S(X)$, be the support of a convexely derived measure $\nu \in \mathcal{MC}_{\rho}$ and let $B=B(x,r)\cap S$. We projectively map $B$ to $\mathbb{R}^n$ and we choose as the center of this projection to be the point $x$. Hence it follows again that
\begin{eqnarray*}
\nu(B)\leq Cr/\rho.
\end{eqnarray*}
\begin{lem}
\label{equi}
Let $\rho>0$. Let 
$\mathcal{K}$ be a compact set of probability measures on $S(X)$ with the following 
property : for every $\nu\in\mathcal{K}$, all $x$ and all $r<\rho$, $\nu(\partial 
B(x,r))=0$.
Then the function $(\nu,x,r)\mapsto
\nu(B(x,r))$ is uniformly continuous on
$\mathcal{K}\times S(X) \times (0,\rho)$. It follows that it is continuous on $\mathcal{MC}^{+}\times S(X) \times[0,\rho)$.
\end{lem}

Proof.

Let $(\nu_i,x_i,r_i)\to (\nu,x,r)$. Let $\{x'_i\}$ (resp $x'$) be the sequence of points (resp the point) on $\mathbb{S}^n$ image of radial projection of the sequence $\{x_i\}$ (resp $x$). Let $\phi_i\in Iso(\mathbb{R}^{n+1})$ be such that $\lim_{i\to \infty}\phi_i=Id$ and for every $i$, $\phi_i(x'_i)=x'$. Such a sequence of isometry acts on $S(X)$ by taking the action on $\mathbb{S}^n$ and projecting to $S(X)$. 
For every $\delta>0$, for big enough $i$ we have
\begin{eqnarray*}
B(x,r-\delta)\subset \phi_i(B(x_i,r_i))\subset B(x,r+\delta).
\end{eqnarray*}
This implies
\begin{eqnarray*}
\nu_i(\phi^{-1}_i(B(x,r-\delta)))<\nu_i(B(x_i,r_i))<\nu_i(\phi^{-1}_i(B(x,r+\delta))).
\end{eqnarray*}
Hence
\begin{eqnarray*}
\limsup \nu_i(B(x_i,r_i))<\lim_{i\to \infty}\phi_{i*}\nu_i(B(x,r+\delta))=\nu(B(x,r+\delta))\\
\liminf \nu_i(B(x_i,r_i))>\lim_{i\to \infty}\phi_{i*}\nu_i(B(x,r-\delta))=\nu(B(x,r-\delta)).
\end{eqnarray*}
Let $\delta \to 0$. As we supposed the norm being smooth, we know that the $\nu(\partial B(x,r))=0$. We can apply the Lemma \ref{radon} and deduce that $\lim_{\delta\to 0}\nu(B(x,r+\delta))=\nu(B(x,r))$. We can apply the Lemma \ref{maj} and the continuity on $\mathcal{MC}^{+}\times S(X) \times[0,\rho)$ is deduced.

\begin{de}[limits of finite convex partitions]
\label{toppartitions}
Let $\Pi$ be a finite convex partition of $S(X)$. We view it as an atomic probability measure $m(\Pi)$ on $\mathcal{MC}$ as follows: for each piece $S$ of $\Pi$, let $\mu_S =\mu_{|S}/\mu(S)$ be the normalized volume of $S$. Then set
\begin{eqnarray*}
m(\Pi)=\sum_{\mathrm{pieces}\,S} \mu(S)\delta_{\mu_S}.
\end{eqnarray*}
We define the {\em space of} (infinite) {\em convex partitions} $\mathcal{CP}$ as the vague closure of the image of the map $m$ in the space $\mathcal{P}(\mathcal{MC})$ of probability measures on the space of convexely derived measures. The subset $\mathcal{CP}^{\leq k}$ of convex partitions of dimension $\leq k$, consists of elements of $\mathcal{CP}$ which are supported on the subset $\mathcal{MC}^{\leq k}$ of convexely derived measures with support of dimension at most $k$.
\end{de}
Note that $\mathcal{CP}$ is compact and $\mathcal{CP}^{\leq k}$ is closed in it. Measures in the support of a convex partition can be thought of as the pieces of the partition.

\begin{lem}[desintegration formula]
\label{desint}
Let $A\subset S(X)$ be a set such that the intersection of $\partial A$ with every $\ell$-dimensional subsphere has vanishing $\ell$-dimensional measure, for all $\ell$, $0<\ell<n$. Let $\Pi\in \mathcal{CP}$. Assume that $\Pi(\mathcal{MC}^{0})=0$. Then
\begin{eqnarray*}
\mu(A)=\int_{\mathcal{MC}}\nu(A)\,d\Pi(\nu).
\end{eqnarray*}
\end{lem}

Proof.

For finite partitions $\Pi_{i}$, equality holds. According to
Lemma \ref{radon}, the function $\nu\mapsto\nu(A)\chi(\nu)$ is continuous
on $\mathcal{MC}^{+}$. Therefore the identity still holds for vague limits of finite partitions. This completes the proof of Lemma.

\subsection{Choice of a Center Map}

In the previous sections we didn't make any particular assumption about the center map. In fact the only property of this map which was used was continuity. In this section we construct a family of center maps which will lead us to the proof of Theorem \ref{usphere}.

\begin{de}[approximate centers of convexely derived measures]\label{centermaps}
Let $\nu \in \mathcal{MC}$, let $r>0$. Consider the function $S(X) \to\mathbb{R}$, $x\mapsto v_{r,\nu} (x)=\nu(B(x,r))$. Let $M_r (\nu)$ be the set of points where $v_{r,\nu}$ achieves its maximum on $\mathrm{support}(\nu)$.

If the support of $\nu$ is $\ell$-dimensional, $\ell<n$, we denote by $M_0(\nu)$ the unique point where the density of $\nu$ achieves its maximum.
\end{de} 

The next Lemma states a semi-continuity property of $M_r$. 

\begin{nott}
When $A_i$, $i\in\mathbb{N}$, are subsets of a topological space, we shall denote by
\begin{eqnarray*}
\lim_{i\to\infty} A_i =\bigcap_{i}\overline{\bigcup_{j\geq i} A_j }.
\end{eqnarray*}
the set of all possible limits of subsequences $x_{i(j)} \in A_{i(j)}$.
\end{nott}

\begin{lem}
\label{semi}
Let $\nu_i$ be convexely derived measures which converge to $\nu \in\mathcal{MC}$. Then, for all $r>0$, 
\begin{eqnarray*}
\lim_{i\to\infty} M_r (\nu_i)\subset M_r (\nu).
\end{eqnarray*}
If follows that
\begin{eqnarray*}
\lim_{i\to\infty} \mathrm{conv.\,hull}(M_r (\nu_i))\subset \mathrm{conv.\,hull}(M_r (\nu)).
\end{eqnarray*}
\end{lem}

Proof.

Let $\nu_i$ tend to $\nu$. Then the support of $\nu_{i}$ Hausdorff
converges to the support of $\nu$. If $\nu\in\mathcal{MC}^{0}$ equals Dirac
measure at $x$, then $M_r (\nu_i)$ automatically converges to $\{x\}=M_r
(\nu)$. Otherwise, $\nu\in\mathcal{MC}^{+}$. Let $x\in \lim_{i\to\infty} M_r (\nu_i)$, i.e. $x=\lim_{i\to\infty}x_i$ for some $x_i \in M_r (\nu_i)$. Pick $y\in\mathrm{support}(\nu)$. Pick a sequence $y_i \in\mathrm{support}(\nu_i)$ converging to $y$. According to Lemma \ref{equi},
\begin{eqnarray*}
v_{r,\nu}(x)=\lim_{i\to\infty} v_{r,\nu_i}(x_i),\quad v_{r,\nu}(y)=\lim_{i\to\infty} v_{r,\nu_i}(y_i).
\end{eqnarray*}
Since $v_{r,\nu_i}(x_i)\geq v_{r,\nu_i}(y_i)$, we get $v_{r,\nu}(x)\geq v_{r,\nu}(y)$, showing that $x\in M_r (\nu)$.

We claim that for arbitrary compact sets $A_i \in S(X)$, 
$$
\lim_{i\to\infty}\mathrm{conv.\,hull}(A_i)\subset\mathrm{conv.\,hull}(\lim_{i\to\infty}A_i). 
$$
Indeed, taking cones, it suffices to check this in $\mathbb{R}^{n+1}$. If $x\in\lim_{i\to\infty}\mathrm{conv.\,hull}(A_i)$, $x=\lim x_i$ with $x_i \in\mathrm{conv.\,hull}(A_i)$, then there exist $n+1$ numbers $t_{i,j}\in[0,1]$ and points $a_{i,j}\in A_i$ such that $\sum_j t_{i,j}=1$, $x_i=\sum_{j}t_{i,j}a_{i,j}$. One can assume that all sequences $i\mapsto t_{i,j}$, $a_{i,j}$ converge to $t_j$, $a_j$. Then $t_j \in[0,1]$, $\sum_j t_j =1$, $a_j \in A=\lim_{i\to\infty}A_i$ and $x=\sum_j t_j a_j \in \mathrm{conv.\,hull}(A)$. This completes the proof of Lemma \ref{semi}.

\medskip

The above semi-continuity property is sufficient to apply Ernest Michael's theory of continuous selections, \cite{Michael}.

\begin{theo}[Michael continuous selection Theorem]\label{michaels}
Let $X$ be paracompact, $Y$ be a Banach Space, and $\mathfrak{S}$ the space of closed convex non-empty subsets of $Y$. Then every lower semi-continuous map $\phi:X\to \mathfrak{S}$ admits a continuous selection.
\end{theo}

Let $N_r(S)$ denotes the convex hull of $M_r(S)$ in $\mathbb{R}^{n+1}$. By definition of the convex sets in $S(X)$, $N_r(S)$ is a closed convex set which does not contain the origin of $\mathbb{R}^{n+1}$. We apply the Theorem \ref{michaels} to the map $S\to N_r(S)$. We obtain a continuous map $S\to D_r(S)$ which never takes the value $0$. We pose $C_r(S)=D_r(S)/\Vert D_r(S)\Vert$ and we obtain the continuous selection on $S(X)$. Hence the following

\begin{de}[centers of open convex sets]\label{centermapsmichael}
Let $r>0$. According to Theorem \ref {michaels}, we can choose a continuous map $C_r :\mathcal{MC}^n \to S(X)$, such that for every $S\in \mathcal{MC}^n$, $C_r (S)$ belongs to $\mathrm{conv.\,hull}(M_r (S))$.
\end{de} 

\subsection{Construction of Partitions Adapted to a Continuous Map}
\begin{de}[partitions adapted to a continuous map]
\label{f}
Let $f:S(X) \to \mathbb{R}^k$ be a continuous map. Let $r\geq 0$. Say a convex partition $\Pi\in\mathcal{CP}$ is {\em $r$-adapted to $f$} if there exists $z\in \mathbb{R}^k$ such that $f^{-1}(z)$ intersects the convex hull of $M_r (\nu)$ for all measures $\nu$ in the support of $\Pi$. Let 
\begin{eqnarray*}
\mathcal{F}_{r}=\{\Pi\in\mathcal{CP}\,|\,\bigcap_{\nu\in\mathrm{support}(\Pi)}f(\mathrm{conv.\,hull}(M_r (\nu)))\not=\emptyset\}
\end{eqnarray*}
denote the set of partitions which are $r$-adapted to $f$.
\end{de}

\begin{prop}
\label{ferme}
For all $r>0$, $\mathcal{F}_r$ is closed in $\mathcal{CP}$. 
\end{prop}

Proof.

If $\lim_{i\to\infty}\Pi_i =\Pi$, $\mathrm{support}(\Pi)\subset\lim_{i\to\infty}\mathrm{support}(\Pi_{i})$, i.e. every piece $\nu$ of $\Pi$ is the limit of a sequence of pieces $\nu_i$ of $\Pi_i$. By assumption, there is a $z_i \in\mathbb{R}^k$ which belongs to all $f(\mathrm{conv.\,hull}(M_r (\nu)))$, $\nu\in\mathrm{support}(\Pi_{i})$. One can assume $z_i$ converges to $z$. Then $z$ belongs to all $f(\mathrm{conv.\,hull}(M_r (\nu)))$, $\nu\in\mathrm{support}(\Pi)$. Indeed, in general, if $g$ is a continuous map and $A_i$ are subsets of a compact space, $g(\lim_{i\to\infty}A_i)=\lim_{i\to\infty}g(A_i)$. So if $\nu=\lim\nu_i$, $\nu_i \in \mathrm{support}(\Pi_{i})$,
\begin{eqnarray*}
z=\lim_{i\to\infty}z_i &\in& \lim_{i\to\infty}f(\mathrm{conv.\,hull}(M_r (\nu_i)))\\
&\subset& f\left(\lim_{i\to\infty}\mathrm{conv.\,hull}(M_r (\nu_i))\right)\\
&\subset& f(\mathrm{conv.\,hull}(M_r (\nu))),
\end{eqnarray*}
thanks to Lemma \ref{semi}.

\begin{cor}
\label{existsfr}
Let $f:S(X) \to \mathbb{R}^k$ be a continuous map. For all $r>0$,
$\mathcal{F}_{r}\cap\mathcal{CP}^{\leq k}$ is non empty.
\end{cor}

Proof.

Theorem 3 states that for every $r>0$, $\mathcal{F}_{r}$ contains uniform atomic measures with arbitrarily many pieces. Theorem 4 produces elements of $\mathcal{F}_{r}$ whose support is contained in arbitrary thin neighborhoods of the compact subset $\mathcal{MC}^{\leq k}$. With Proposition \ref{ferme}, this gives elements in $\mathcal{F}_{r}\cap\mathcal{CP}^{\leq k}$.
\subsection{Convergence of $M_r(\nu)$ as $r$ tends to $0$}

\begin{lem}
\label{convargmax}
Let $\ell<n$. For every $\ell$-dimensional convexely derived measure $\nu$,
\begin{eqnarray*}
\lim_{r\to 0}d_{H}(M_r (\nu),M_0 (\nu))=0.
\end{eqnarray*}
\end{lem}

Proof.

We prove the Lemma by contradiction. Otherwise, we get a $\delta>0$ and a
sequence of radii $r_i$ tending to $0$ such that $d_{H}(M_{r_i}(\nu),M_0 (\nu))\geq\delta$. Pick a point $x_{i} \in S$ where $v_{r_i,\nu}$ achieves its maximum and such that $d(x_{i},M_0 (\nu))\geq\delta$. Up to extracting a subsequence, we can assume that $x_{i}$ converges to $x\in S$. Then $v_{r_i,\nu}(x_i)/\alpha_k r_i^k$ converges to $\phi_{\nu}(x)$. For every $y\in S$, $v_{r_i,\nu}(y)\leq v_{r_i,\nu}(x)$ and $v_{r_i,\nu}(y)/\alpha_k r_i^k$ converges to $\phi_{\nu}(y)$. Therefore $\phi_{\nu}(y)\leq\phi_{\nu}(x)$. This shows that $\{x\}=M_0 (\nu)$, contradiction. 

\bigskip

A stronger statement (Corollary \ref{unifgood}) will be given after the following technical lemmas.

\begin{lem}
\label{boundeddensity}
Let $\nu$ be a convexely derived measure on $S(X)$ whose support is a $k$-dimensional convex set $S$. Write $d\nu=\phi\,d\mu_k$. Then
\begin{eqnarray*}
\max_{S}\phi\leq \frac{2^{n+1}}{\mu_{k}(S)}.
\end{eqnarray*}
\end{lem}

Proof.

Replace $S$ with $C=co(S)\subset\mathbb{R}^{n+1}$, and $\phi$ by its $n-k$-homogeneous extension. Then $\phi^{1/(n-k)}$ is concave. Assume $\phi$ achieves its maximum at $x\in C$. Translate $C$ so that $x=0$. On $\frac{1}{2}C$, $\phi^{1/(n-k)}\geq\frac{1}{2}\phi^{1/(n-k)}(x)$, thus
\begin{eqnarray*}
1=\nu(S)&\geq&\int_{\frac{1}{2}C}\phi\,dvol_{k+1}\\
&\geq&\frac{1}{2^{n-k}}\phi(x)vol_{k+1}(\frac{1}{2}C)\\
&=&\frac{1}{2^{n+1}}\phi(x)vol_{k+1}(C)\\
&=&\frac{1}{2^{n+1}}\phi(x)\mu_{k}(S).
\end{eqnarray*}

\begin{lem}
\label{equiconc}
Let $S$, $S_i$ be full compact convex subsets of $\mathbb{R}^{n}$ such that $S_i$ Hausdorff-converges to $S$. Let $\phi_i:S_i \to [0,1]$ be concave functions. Then there exists a concave function $\phi:S\to[0,1]$ and a subsequence with the following properties.
\begin{itemize}
  \item On every compact subset of the interior of $S$, $\phi_i$ converges uniformly to $\phi$.
  \item For all $x\in\partial S$ and all sequences $x_i \in S_i$ converging to $x$, 
\begin{eqnarray*}
\limsup_{i\to\infty}\phi_i (x_i) \leq \phi(x).
\end{eqnarray*}
\end{itemize}
\end{lem}

Proof.

In general, bounded concave functions $f$ on compact convex sets $\Sigma$ are locally Lipschitz, 
\begin{center}
{\em for $x\in \Sigma$ with $d(x,\partial \Sigma)=r$, and all $y\in \Sigma$, $\displaystyle |f(x)-f(y)|\leq\frac{1}{r}d(x,y)$}. 
\end{center}
Indeed, let $[x',y']$ be the intersection of $\Sigma$ with the line through $x$ and $y$, with $x'$, $x$, $y'$ and $y'$ sitting along the line in this order. Let $\ell$ be the affine function on $[x',y']$ such that $\ell(x')=f(x')$ and $\ell(x)=f(x)$. Then $f(y)\leq \ell(y)$, thus $f(y)-f(x)\leq\frac{1}{d(x',x)}|f(x)-f(x')|d(x,y)\leq\frac{1}{r}d(x,y)$. Also, let $\ell'$ be the affine function on $[x',y']$ such that $\ell'(x)=f(x)$ and $\ell'(y')=f(y')$. Then $f(y)\geq \ell'(y)$, thus $f(y)-f(x)\geq-\frac{1}{d(x,y')}|f(x)-f(y')|d(x,y)\geq-\frac{1}{r}d(x,y)$. 

This shows that on every compact subset of the interior of $S$, the sequence $f_j$ is equicontinuous, so a subsequence can be found which converges uniformly on all such compact sets to a continuous function $\phi$. Of course, $\phi$ is concave and bounded, so it extends continuously to $\partial S$. Let $x\in\partial S$ and $x_i \in S_i$ converge to $x$. Pick an interior point $x_0$ of $S$ and a second interior point $x'\not=x_0$ such that $x_0$ lies on the segment $[x',x]$. Pick $x'_i$ on the line passing through $x_0$ and $x_i$ and converging to $x'$. The Lipschitz estimate for $\phi_i$ reads
\begin{eqnarray*}
\phi_i (x_i)-\phi_i (x_0)\leq \frac{d(x_0 ,x_i)}{d(x_0 ,x'_i)}|\phi_i (x'_i)-\phi_i (x_0)|.
\end{eqnarray*}
Letting $i$ tend to infinity yields
\begin{eqnarray*}
\limsup\phi_i (x_i)\leq\phi(x_0) + \frac{d(x_0 ,x)}{d(x_0 ,x')}|\phi(x')-\phi (x_0)|.
\end{eqnarray*}
Letting $x_0$ and $x'$ tend to $x$ (while keeping $x'$, $x_0$ and $x$ aligned and $ \frac{d(x_0 ,x)}{d(x_0 ,x')}$ bounded) gives $\limsup\phi_i (x_i)\leq\phi(x)$.

\begin{lem}
\label{good}
For each $k<n$, the restriction of $(\nu,r)\mapsto d_H (M_{r}(\nu),M_0 (\nu))$ to $\mathbb{R}_+ \times \mathcal{MC}^{k}$ tends to $0$ along $\{0\} \times \mathcal{MC}^{k}$, i.e. for all $\nu\in\mathcal{MC}^{k}$,
\begin{eqnarray*}
\lim_{r\to 0,\,\nu'\to\nu,\, \nu'\in \mathcal{MC}^{k}}d_H (M_{r}(\nu),M_0 (\nu))=0.
\end{eqnarray*}
\end{lem}

Proof.

Let $\nu\in \mathcal{MC}^{k}$. Let $\nu_i$ be a sequence of $k$-dimensional convexely derived measures which converges to  $\nu$ and $r_i$ be  positive numbers tending to $0$. For every $i$, we project the support of $\nu_i$ into the $k$-sphere which contains the support of $\nu$ (if intrinsically this poses problem, one can always think of the cones over the support of these measure and do all projections in $\mathbb{R}^{n+1}$). In other words, one can assume that all $\nu_i$ have support $S_i$ in the same $k$-sphere. Of course, $S_i$ Hausdorff-converges to the support $S$ of $\nu$. Let $\phi_i$ denote the density of $\nu_i$ with respect to $k$-dimensional conical measure. Since $\mu_{k}(S_i)$ does not tend to $0$, $\phi_i$ are uniformly bounded, by Lemma \ref{boundeddensity}. Furthermore, on any compact convex subset $K$ of the relative interior of $S$, the $\phi_i$ are equicontinuous (this follows by the cone construction from Lemma \ref{equiconc}). Therefore one can assume that $\phi_i$ converge uniformly on compact subsets of the relative interior of $S$. Since for all $r'>0$, $v_{r',\nu_i}$ converges to $v_{r',\nu}$, the limit must be equal to the density $\phi$ of $\nu$. From Lemma \ref{equiconc}, one can assert that at boundary points $x\in\partial S$, for every sequence $x_i \in S_i$ converging to $x$, $\limsup \phi_i (x_i)\leq \phi(x)$.

We repeat the argument of Lemma \ref{convargmax}. If $M_{r_i}(\nu_i)$ does not converge to $M_0 (\nu)$, some sequence $x_i \in M_{r_i}(\nu_i)$ satisfies $d(x_{i},M_0 (\nu))\geq\delta$ for some $\delta>0$. Up to extracting a subsequence, we can assume that $x_{i}$ converges to $x\in S$. If $x\notin\partial S$, then $v_{r_i,\nu}(x_i)/\alpha_k r_i^k$ converges to $\phi(x)$. If $x\in\partial S$, $\limsup v_{r_i,\nu}(x_i)/\alpha_k r_i^k \leq \phi(x)$. For every $y\in S\setminus\partial S$, $v_{r_i,\nu}(y)\leq v_{r_i,\nu}(x)$ and $v_{r_i,\nu}(y)/\alpha_k r_i^k$ converges to $\phi(y)$. Therefore $\phi(y)\leq\phi(x)$. Since $S\setminus\partial S$ is dense in $S$, this holds for all $y\in S$, thus $\phi$ achieves its maximum at $x$, i.e. $\{x\}=M_0 (\nu)$, contradiction.

\begin{cor}
\label{unifgood}
On any compact subset of $\mathcal{MC}^{k}$, the functions 
$$
\nu\mapsto d_H (M_{r}(\nu),M_0 (\nu))
$$ 
converge uniformly to $0$ as $r$ tends to $0$.
\end{cor}

\begin{prop}
\label{goodcase}
Assume $f:S(X) \to \mathbb{R}^k$ is a generic smooth map. Let $r_i$ tend to $0$ and let $\Pi_i \in\mathcal{CP}^{\leq k}\cap\mathcal{F}_{r_i}$ be convex partitions of dimension $\leq k$,  $r_i$-adapted to $f$. Then, for all $\varepsilon>0$,
\begin{eqnarray*}
\max_{z\in\mathbb{R}^k}\mu(f^{-1}(z)+\varepsilon)\geq w(\varepsilon)\limsup_{i\to\infty}\Pi_{i}(\mathcal{MC}^{k}).
\end{eqnarray*}
Where 
\begin{eqnarray*}
w(\varepsilon)= \frac{1}{1+(1-2\delta(\frac{\varepsilon}{2}))^{n-k}(k+1)^{k+1}\frac{F(k,\frac{\varepsilon}{2})}{G(k,\frac{\varepsilon}{2})}}.
\end{eqnarray*}
And where the functions $F(.,.)$ and $G(.,.)$ were defined previously.
\end{prop}

Proof.

By assumption, for each $i$, there exists $z_i \in \mathbb{R}^k$ such that for all $\mu\in\mathrm{support}(\Pi_i)$, there exists $x_{i,\nu}\in \mathrm{conv.\,hull}(M_{r_i}(\nu))$ such that $f(x_{i,\nu})=z_i$. Let $\mathcal{K}\subset \mathcal{MC}^{k}$ be a compact set. According to Corollary \ref{unifgood} and Lemma \ref{equi}, for all $\varepsilon>0$,
\begin{eqnarray*}
\delta_i :=\sup_{\nu\in\mathcal{K}}|\nu(B(x_{i,\nu},\varepsilon))-\nu(B(M_{0}(\nu),\varepsilon))|
\end{eqnarray*}
tends to $0$. Considerations in previous sections show that for every $k$-dimensional convexely derived measure $\nu$, 
$$
\nu(B(M_{0}(\nu),\varepsilon))\geq w(\varepsilon).
$$
For a generic smooth map $f$, the intersection of $f^{-1}(z_i)+\varepsilon$ with $k$-dimensional convex sets has vanishing $k$-dimensional measure, so the desintegration formula applies, and
\begin{eqnarray*}
\mu(f^{-1}(z_i)+\varepsilon)
&\geq&\int_{\mathcal{MC}^{+}}\nu(f^{-1}(z_i)+\varepsilon)\,d\Pi_i (\nu)\\
&\geq&\int_{\mathcal{K}}\nu(B(x_{i,\nu},\varepsilon))\,d\Pi_i (\nu)\\
&\geq&\Pi_i (\mathcal{K})w(\varepsilon)-\delta_i .
\end{eqnarray*}
Taking the supremum over all compact subsets of $\mathcal{MC}^{k}$ and then a limit as $i$ tends to infinity yields the announced inequality.

\medskip

\subsection{End of the Proof of Theorem \ref{usphere}}
There remains to show that convex partitions in $\mathcal{CP}^{\leq k}\cap\mathcal{F}_{r}$, $r$ small, put most of their weight on $k$-dimensional pieces. This will be proven indirectly. Pieces of dimension $<k$ may exist, but they provide a lower bound on $\mu(f^{-1}(z)+r)$ which is so large, that they must have small weight.
We shall need a weak concavity property of $v_{\mu,r}$, which in turn relies on the corresponding Euclidean statement.

\begin{lem}
\label{gpl}
Let $S \subset \mathbb{R}^n$ be an open convex set, $\phi$ an $m$-concave function defined on $S$. Let $\mu=\phi dvol_n$. Then the map $x\mapsto\mu(B(x,r)\cap S)$ is $m+n$-concave on $S$.
\end{lem}

Proof.

We use the following estimate (Generalized Prekopa-Leindler inequality), which can be found in \cite{led}. For $\alpha\in[-\infty,+\infty]$ and $\theta\in[0,1]$, the $\alpha$-mean of two nonnegative numbers $a$ and $b$ with weight $\theta$ is
\begin{eqnarray*}
M_{\alpha}^{(\theta)}(a,b)=(\theta a^\alpha +(1-\theta)b^\alpha)^{1/\alpha}.
\end{eqnarray*}
Let $-\frac{1}{n}\leq\alpha\leq+\infty$, $\theta\in[0,1]$, $u$, $v$, $w$ nonnegative measurable functions on $\mathbb{R}^n$ such that for all $x$, $y\in\mathbb{R}^n$,
\begin{eqnarray*}
w(\theta x+(1-\theta)y)\geq M_{\alpha}^{(\theta)}(u(x),v(y)).
\end{eqnarray*}
Let $\beta=\frac{\alpha}{1+\alpha n}$. Then
\begin{eqnarray*}
\int w\geq M_{\beta}^{(\theta)}(\int u,\int v).
\end{eqnarray*}
We apply this to restrictions of $\phi$ to balls, $u=1_{B(x,r)}\phi$, $v=1_{B(y,r)}\phi$, $w=1_{B(\theta x+(1-\theta)y,r)}\phi$. By $m$-convexity of $\phi$, the assumptions of the generalized Prekopa-Leindler inequality are satisfied with $\alpha=1/m$. Then for $\beta=\frac{1}{m+n}$,
\begin{eqnarray*}
\mu(B(\theta x+(1-\theta)y),r))\geq M_{\beta}^{(\theta)}(\mu(B(x,r)),\mu(B(y,r))),
\end{eqnarray*}
which means
\begin{eqnarray*}
\mu(B(\theta x+(1-\theta)y),r))^{\frac{1}{m+n}}\geq \theta\mu(B(x,r))^{\frac{1}{m+n}}+(1-\theta)\mu(B(y,r))^{\frac{1}{m+n}}.
\end{eqnarray*}
\begin{lem}
\label{wconcave}
The functions $v_{\nu,r}$ are weakly concave on $S(X)$. In other words, there exists
a constant $c=c(n)>0$ such that for every convexely derived measure $\nu$
and every sufficiently small $r>0$, if $K\subset\mathrm{support}(\nu)$, then
\begin{eqnarray*}
\min_{\mathrm{conv}(K)}v_{\nu,\frac{r}{c}}\geq c\,\min_{K}v_{\nu,r}.
\end{eqnarray*}
\end{lem}

Proof.

Since a half-sphere is projectively equivalent with Euclidean space, it suffices to prove weak concavity when $K$ consists of 2 points.

Let $\nu$ be a $k$-dimensional convexely derived measure on $S(X)$. Denote its density by $\phi$, a weak $(n-k)$-concave function on the support $S$ of $\nu$. Let $\Phi$ denote the $(n-k)$-homogeneous extension of $\phi$ to the cone on $S$. This is $(n-k)$-concave. Fix a point $x_0 \in S(X)$, let $\mathbb{R}^n$ denote the tangent space (cone) of $S(X)$ at $x_0$. Denote by $\phi'$ the restriction of $\Phi$ to $\mathbb{R}^n$, and $\nu'$ the measure with density $\phi'$. Lemma \ref{gpl} implies that $x'\mapsto \mu(B(x',r))$ is $(2n-k)$-concave. This implies that for every $x'$, $y'\in\mathbb{R}^n$ and $z'$ belonging to the middle third of the line segment $[x',y']$,
\begin{eqnarray*}
\nu'(B(z',r))\geq \frac{1}{3^{2n-k}}\max\{\nu'(B(x',r)),\nu'(B(y',r))\}.
\end{eqnarray*}

The radial projection from a neighborhood $V\subset S(X)$ of $x_0$ to $\mathbb{R}^n$ is nearly isometric and nearly maps $\phi'$ to $\phi$. Thus there exists a constant $c_1 >0$ such that if $x$, $y\in V$ and $z$ belongs to the middle third of the segment $[x,y]$,
\begin{eqnarray*}
\nu(B(z,\frac{r}{c_1}))\geq c_1 \,\max\{\nu(B(x,r)),\nu(B(y,r))\}.
\end{eqnarray*}
Covering long segments $[x,y]$ with $N$ neighborhoods like $V$ ($N$ can be bounded independantly of $n$) provides a constant $c>0$ such that for all $z\in[x,y]$ which is not too close to the endpoints,
\begin{eqnarray*}
\nu(B(z,\frac{r}{c_1^N}))\geq c_1^N \,\max\{\nu(B(x,r)),\nu(B(y,r))\}.
\end{eqnarray*}
In particular, for $c=c_1^N$,
\begin{eqnarray*}
\nu(B(z,\frac{r}{c}))\geq c \,\min\{\nu(B(x,r)),\nu(B(y,r))\}.
\end{eqnarray*}

\begin{prop}
\label{badcase}
There exists a constant $c=c(n)>0$ such that if  $f:S(X)
\to\mathbb{R}^k$ is smooth and generic and $\Pi$ belongs to $\mathcal{F}_{r}\cap\mathcal{CP}^{\leq k}$ for some small enough $r>0$, then,
\begin{eqnarray*}
\max_{z\in\mathbb{R}^k}\mu(f^{-1}(z)+\frac{r}{c})\geq c\,\sum_{\ell=0}^{k}w_l(r)\Pi(\mathcal{MC}^{\ell}).
\end{eqnarray*}
Where $w_l(r)$ is equal to $w(r)$ in codimension $l$.
\end{prop}

Proof.

By assumption, there exists $z\in\mathbb{R}^{k}$ such that for every measure $\nu$ in the support of $\Pi$, there exists $x\in \mathrm{conv.\,hull}(M_r (\nu))$ such that $f(x)=z$. If the support of $\nu$ is $\ell$-dimensional, Lemma \ref{semi} and our previous computations
\begin{eqnarray*}
\nu(f^{-1}(z)+\frac{r}{c})&\geq& \nu(B(x,\frac{r}{c}))\\
&=&v_{\mu,\frac{r}{c}}(x)\\
&\geq& c\,\min_{M_r (\nu)}v_{\nu,r}\\
&=&c\,\max_{\mathrm{support}(\nu)}v_{\nu,r}\\
&\geq& c\,v_{\nu,r}(M_0 (\nu)) \\
&=&c\,\nu(B(M_0 (\nu),r))\\
&\geq& c\,w_l(\rho). 
\end{eqnarray*}
Again, for generic smooth $f$, one can integrate this with respect to
$\Pi$. 
\begin{eqnarray*}
\mu(f^{-1}(z)+r)&=&\int_{\mathcal{MC}}\nu(f^{-1}(z)+r)\,d\Pi(\nu)\\
&\geq& c\,\sum_{\ell=0}^{k}w_l(\rho)\Pi(\mathcal{MC}^{\ell}).
\end{eqnarray*}
\begin{lem}
For every $l<k$, we have
\begin{eqnarray*}
\lim_{r\to 0}w_l(r)/w_k(r)=\infty
\end{eqnarray*}
\end{lem}

Proof.

Simple observation shows that for every $m \in \mathbb{N}$, $\lim_{r\to 0}G(m,r)\to 0$, and $\lim_{r\to 0}F(m,r)=1$.
Simple calculation leads to
\begin{eqnarray*}
w_l(r)/w_k(r)=\frac{1+(1-2\delta(r/2))^{n-k}\frac{F(k,r/2)}{G(k,r/2)}(k+1)^{k+1}}{1+(1-2\delta(r/2))^{n-l}\frac{F(l,r/2)}{G(l,r/2)}(l+1)^{l+1}}\backsim_{r \to 0} C\frac{G(l,r)}{G(k,r)}
\end{eqnarray*}
And by the well known asymptotic behvior of the function $G(m,r)$ we have
\begin{eqnarray*}
\frac{G(l,r)}{G(k,r)} \backsim_{r \to 0} r^{l-k}.
\end{eqnarray*}
Hence the proof of the Lemma follows.

\textbf{Proof of Theorem \ref{usphere}}.

\begin{prop}
\label{5}
Let $\varepsilon>0$. Let $f:S(X) \to\mathbb{R}^k$ be a continuous map. Then
\begin{eqnarray*}
\max_{z\in\mathbb{R}^k}\mu(f^{-1}(z)+\varepsilon)\geq w(\varepsilon).
\end{eqnarray*}
\end{prop}

Proof.

Assume first that $f$ is smooth and generic. Then there exists a constant $W$ such that for all sufficiently small $r$,
\begin{eqnarray*}
\max_{z\in\mathbb{R}^k}\mu(f^{-1}(z)+r)\leq Wr^{k}.
\end{eqnarray*}
For every $r>0$, there exists a convex partition $\Pi_r \in
\mathcal{CP}^{\leq k}\cap\mathcal{F}_r$ which is $r$-adapted to $f$
(Corollary \ref{existsfr}). Proposition \ref{badcase} yields
\begin{eqnarray*}
\sum_{\ell=0}^{k}w_l(r)\Pi_r (\mathcal{MC}^{\ell})\leq \frac{1}{c}\,\max_{z\in\mathbb{R}^k}\mu(f^{-1}(z)+\frac{r}{c})\leq \frac{W}{c}(\frac{r}{c})^k.
\end{eqnarray*}
As $r$ tends to $0$, this implies that for all $\ell<k$ (including $\ell=0$), $\Pi_r (\mathcal{MC}^{\ell})$ tends to $0$, and thus $\Pi_r (\mathcal{MC}^{k})$ tends to $1$. Letting $r$ tend to $0$ in Proposition \ref{goodcase} then shows that 
\begin{eqnarray*}
\max_{z\in\mathbb{R}^k}\mu(f^{-1}(z)+\varepsilon)\geq w(\varepsilon).
\end{eqnarray*}

Every continuous map $f:S(X) \to\mathbb{R}^k$ is a uniform limit of smooth generic maps. Hausdorff semi-continuity of $X\mapsto \mu(X+\varepsilon)$ then extends the result to all continuous maps. Indeed, let the continuous map $f: S(X) \to \mathbb{R}^k$ of theorem \ref{usphere} be fixed. Let $g_j :S(X) \to \mathbb{R}^k$ be a sequence of $C^{\infty}$ maps such that $\delta_{j}=\Vert g_j -f \Vert_{C^0}$ tends to $0$. For every $j$, there exists a $z_j \in \mathbb{R}^k$ such that $\mu(g_j^{-1}(z_j)+\varepsilon) \geq w(\varepsilon)$. We know that for every $j$, $g_{j}^{-1}(z_j) \subseteq f^{-1}(B(z_j,\delta_{j}))$. Then
$$\mu (f^{-1}(B(z_j,\delta_{j}))+\varepsilon) \geq \mu (g_{j}^{-1}(z_j)+\varepsilon) \geq w(\varepsilon).$$
Up to extracting a subsequence, we can assume that $\{z_j\}$ converges to a point $z$. There exists a decreasing sequence $\varepsilon_j \to 0$ such that for every $j$, $\vert z-z_j \vert \leq \varepsilon_j$. Then
$$f^{-1}(B(z_j,\delta_{j}))+\varepsilon \subseteq f^{-1}(B(z,\delta_{j}+\varepsilon_j))+\varepsilon,$$
thus for all $j$
$$\mu(f^{-1}(B(z,\delta_{j}+\varepsilon_j)+\varepsilon) \geq w(\varepsilon),$$
and by Fatou Lemma
$$\mu(\underset{j}{\bigcap}f^{-1}(B(z,\delta_{j}+\varepsilon_j))+\varepsilon) \geq w(\varepsilon).$$
If for all $j$, $x \in f^{-1}(B(z,\delta_{j}+\varepsilon_j))+\varepsilon$, then there exists $y_j$ such that $d(x,y_j) \leq \varepsilon$ and $f(y_j) \in B(z,\delta_{j}+\varepsilon_j)$. We choose a subsequence $y_k$ which converges to $y$. By construction, $d(x,y) \leq \varepsilon$, $f(y)=z$ thus $x \in f^{-1}(z)+\varepsilon$. Hence
\begin{eqnarray*}
\underset{j}{\bigcap}f^{-1}(B(z,\delta_{j}+\varepsilon_j))+\varepsilon)\subset f^{-1}(z)+\varepsilon,
\end{eqnarray*}
and
\begin{eqnarray*}
\mu(f^{-1}(z)+\varepsilon) \geq w(\varepsilon).
\end{eqnarray*}

\section{Why all these complications?}
Remember the following
\begin{theo}[Gromov 2003]
\label{1.1}
Let $f : \mathbb{S}^n \to \mathbb{R}^k$ be a continuous map from the canonical unit $n$-sphere to a Euclidean space of dimension $k$ where $k \leq n$. There exists a point $z \in \mathbb{R}^k$ such that the $n$-spherical volume of the $\varepsilon$- tubular neighborhood of $f^{-1}(z)$, denoted by $f^{-1}(z)+ \varepsilon$ satisfies, for every $\varepsilon > 0$,
$$ vol_n (f^{-1}(z)+ \varepsilon) \geq vol_n (S^{n-k}+\varepsilon).$$
Here $\mathbb{S}^{n-k}$ is the $(n-k)$-equatorial sphere of $\mathbb{S}^n$.
\end{theo}
Several times during the last sections, we used the radial projection between the canonical sphere and the unit sphere $S(X)$. One (including myself) could ask why bothering with all we did and not just radially projecting the result of Theorem \ref{1.1} on $S(X)$. Indeed, this gives another lower bound for the waist of $S(X)$ as we will show in the next

\begin{prop}
Let $X$ be a uniformly convex normed space of finite dimension $n+1$. Let $S(X)$ be the unit sphere of $X$, for which the distance is induced from the norm of $X$. The measure defined on $S(X)$ is the conical probability measure. So a lower bound for the waist of $S(X)$ relative to $\mathbb{R}^{k}$ is given by
$$w_2(\varepsilon)=(n+1)^{-n-1}\frac{vol(\mathbb{S}^{n-k}+\frac{\varepsilon}{n+1})}{vol(\mathbb{S}^n)}.$$
\end{prop}

\emph{Proof of the Proposition}

Let $pr$ be the radial projection of $\mathbb{S}^n$ to $S(X)$. We apply theorem \ref{1.1} to the map $g=pr^{-1}\circ f$. Hence there exists a fiber $X$ such that for every $\varepsilon>0$ 
\begin{eqnarray*}
vol(X+\varepsilon)\geq vol(\mathbb{S}^{n-k}+\varepsilon)
\end{eqnarray*}
We radially project $X+\varepsilon$ to $S(X)$. We have
\begin{eqnarray*}
pr(X+\varepsilon)\subset pr(X)+(n+1)\varepsilon
\end{eqnarray*}
Hence
\begin{eqnarray*}
\mu(pr(X)+\varepsilon)&\geq& \mu(pr(X+\frac{\varepsilon}{n+1}))\\
&\geq&(n+1)^{-n-1}\frac{vol(X+\frac{\varepsilon}{n+1})}{vol(\mathbb{S}^n)}\\
&\geq&(n+1)^{-n-1}\frac{vol(\mathbb{S}^{n-k}+\frac{\varepsilon}{n+1})}{vol(\mathbb{S}^n)}.
\end{eqnarray*}
And the proposition is proved.

\bigskip

We see that a brutal application of Gromov's theorem gives a lower bound for the waist of the unit sphere of a uniformly convex normed space, $S(X)$. But comparing $w_1(\varepsilon)$ and $w_2(\varepsilon)$, we can see that the lower bound $w_1(\varepsilon)$ has a much better dependence on the variable $n$, even if the dependence on the variable $k$ is very bad. 

For example, if $k$ is fixed and $n$ tends to infinity, $w_2(\varepsilon)$ tends (exponentially fast) to $0$ while for this case, the lower bound $w_1(\varepsilon)$ tends to $1$. One can hope to have a better dependence on the variable $k$ by knowing the best degree of dilation of the radial projection of $\mathbb{S}^n \to S(X)$. Here we gave a trivial bound for the degree of dilation, not taking into account uniform convexity.

\section{Comparison with Gromov-Milman's isoperimetric inequality}

We need to compare the result of Theorem \ref{usphere} for $k=1$ with Gromov-Milman's isoperimetric inequality which we recall here. This inequality was proved first by Gromov-Millman in \cite{gromil}. The proof was completed later on by S. Alesker in \cite{ale} (S. Sodin had the kindness of referring Alesker's paper to the author). There is a very short and easy proof given by J. Arias-de-Reyna, K. Ball and R. Villa in \cite{ball}.

\begin{theo} \label{grmil}
Let $S(X)$ be a uniformly convex unit sphere with modulus $\delta$. For every Borel set $A \subset S(X)$ such that $\mu(A)\geq \frac{1}{2}$ and for every $\varepsilon>0$ we have
\begin{eqnarray*}
\mu(A+\varepsilon)\geq 1-e^{-a(\varepsilon)n},
\end{eqnarray*}
where $a(\varepsilon)=\delta(\frac{\varepsilon}{8}-\theta_n)$ and where $\theta_n=1-(\frac{1}{2})^{1/(n-1)}$.
\end{theo}

Our Theorem \ref{usphere}, in case $k=1$, implies a similar isoperimetric inequality.

We need the following proposition which relates isoperimetry and $1$-waist.

\begin{prop} \label{waistiso}
$1$-waist $\Rightarrow$ Isoperimetry : For every open subset $A\subset S(X)$ and for all $\varepsilon>0$ we have
\begin{eqnarray*}
\max\{\mu(A+\varepsilon),\mu(A^c+\varepsilon)\}\geq w(\varepsilon).
\end{eqnarray*}
\end{prop}
For the proof, see \cite{memphd} where we prove this Proposition in a more general context. 

Proposition \ref{waistiso} is far from optimal for small $\varepsilon$ and fixed $n$. On the other hand, let $\varepsilon$ be fixed and let $n \to \infty$. In this regime, our main theorem \ref{usphere} combined with Proposition \ref{waistiso} yields
\begin{eqnarray*}
\max\{\mu(A+\varepsilon),\mu(A^{c}+\varepsilon)\}\geq 1-e^{-b(\varepsilon)n-c(\varepsilon)},
\end{eqnarray*}
where $b(\varepsilon)=2\delta(\frac{\varepsilon}{2})$ and $c(\varepsilon)$ has an ugly expression. Since, $b>a$, our theorem \ref{usphere} gives a better estimate.

\bibliographystyle{plain}
\bibliography{bibusphere}

\end{document}